\documentclass[12pt]{article}
\usepackage[left=2.5cm,right=2.5cm,top=2.5cm,bottom=2.5cm,a4paper]{geometry}

\usepackage[a4paper]{geometry}
\usepackage{graphicx}
\usepackage{microtype}
\usepackage{siunitx}
\usepackage{booktabs}
\usepackage{graphics}
\usepackage{graphicx}
\usepackage{epsfig}
\usepackage{amsmath,amsfonts,amssymb,amsthm}
\usepackage[hidelinks]{hyperref}
\usepackage{cleveref}
\usepackage{listings}
\usepackage{paralist}
\usepackage{sectsty}
\usepackage{datetime}
\usepackage{algorithmic,algorithm}
\usepackage[X2,T1]{fontenc} 

\usepackage{authblk}

\theoremstyle{definition}
\newtheorem{definition}{Definition}[section]

\theoremstyle{plain}
\newtheorem{theorem}[definition]{Theorem}

\newtheorem{lemma}[definition]{Lemma}

\title{Two Quintic Diophantine Equations and Fifth Powers in Fibonacci Sums at Distance Five}
\author[1]{Kyunghwan Song\thanks{khsong@jejunu.ac.kr}}
\affil[1]{Department of Mathematics, Jeju National University, 102 Jejudaehakro Jeju, 63243, Republic of Korea}
\date{}

\begin{document}
\maketitle
\begin{abstract}
Let $F_n$ denote the $n$th Fibonacci number. We study the Diophantine
 equation
$$
 F_{n+5}+F_n=y^5
$$
in positive integers $n$ and $y$. We prove that this equation has no
solutions. The proof uses a recurrence identity to reduce the problem
to the two equations
$$
 X^2-5Z^5=44,\qquad X^2+44=5Z^5.
$$
We show that the first equation has exactly the solutions
$(X,Z)=(\pm7,1)$ for $Z>0$, and that the second has no such solutions.
To obtain these results, we use ideal factorization in
$\mathbb Q(\sqrt{11})$ and $\mathbb Q(\sqrt{-11})$ and reduce the
problem to six quintic Thue equations. Their complete solution sets
are computed with the Magma Online Calculator and independently
verified with PARI/GP.
\end{abstract}

\section{Introduction}\label{sec:intro}

Let $F_n$ be the Fibonacci sequence defined by
$$
F_0=0,\qquad F_1=1,\qquad F_{n+2}=F_{n+1}+F_n\quad(n\ge0).
$$
The problem of finding perfect powers in Fibonacci numbers has been
studied by several authors. Bugeaud, Mignotte, and Siksek
\cite{BMS2006} determined all perfect powers in the Fibonacci and
Lucas sequences. Their proof uses methods from algebraic number
theory and the theory of Diophantine equations.

A related problem concerns perfect powers that can be written as the
sum of two Fibonacci numbers. Luca and Patel \cite{LucaPatel2018}
studied
\begin{equation}\label{eq:general}
F_m+F_n=y^p,
\end{equation}
where $m,n,y$ are positive integers and $p\ge2$. In particular, they
obtained a classification when $m$ and $n$ have the same parity.
The work of Florian Luca on perfect powers and Fibonacci sums is
particularly relevant to the present paper. In particular, the
joint work of Luca and Patel provides a natural starting point
for studying fixed differences between Fibonacci indices.
Further results concerning sums of Fibonacci numbers and perfect
powers were given by Ziegler \cite{Ziegler2023}. These studies
motivate the investigation of particular relations between the two
indices.

In this paper, we give a complete determination of perfect fifth
powers in sums of two Fibonacci numbers whose indices differ by five.
More precisely, we consider
\begin{equation}\label{eq:main}
F_{n+5}+F_n=y^5,\qquad n,y\ge1.
\end{equation}
Since $n$ and $n+5$ have different parity, the same-parity result of
Luca and Patel does not directly apply to \eqref{eq:main}.
In this sense, our fixed-gap problem is directly connected to
the Diophantine questions investigated by Florian Luca and his
collaborators. We consider a complementary parity setting,
without claiming that the general opposite-parity problem is resolved. We shall
prove that \eqref{eq:main} has no solutions in positive integers.

Put
$$
S_n=F_{n+5}+F_n.
$$
The sequence $(S_n)$ satisfies the Fibonacci recurrence and has
initial values $S_0=5$ and $S_1=9$. The identity
$$
(2S_{n+1}-S_n)^2-5S_n^2=44(-1)^n
$$
reduces \eqref{eq:main} to the following two Diophantine equations:
\begin{equation}\label{eq:pair}
X^2-5Z^5=44,\qquad X^2+44=5Z^5.
\end{equation}
We determine the solutions of both equations for $Z>0$. This also
gives the main result for the Fibonacci sequence.

The two equations in \eqref{eq:pair} are treated by related methods.
For the first equation, we work in $\mathbb Q(\sqrt{11})$ and obtain
five quintic Thue equations. For the second, we work in
$\mathbb Q(\sqrt{-11})$ and obtain one quintic Thue equation. The
reductions use the class number and unit group of the two quadratic
fields. We solve the resulting Thue equations with the Magma Online
Calculator and check the complete solution sets independently with
PARI/GP. The use of computational Thue solvers in Diophantine
problems has precedents, for example, in \cite{McDaniel2011,BMS2006II}.
The computer calculations are used only after the theoretical
reductions have been established.

We state our main results as follows.
\begin{theorem}\label{thm:diophantine}
For integers $X,Z$ with $Z>0$, the equation
$X^2-5Z^5=44$ has exactly the solutions
$(X,Z)=(7,1)$ and $(-7,1)$. The equation $X^2+44=5Z^5$
has no integer solutions with $Z>0$.
\end{theorem}
\begin{theorem}\label{thm:intro}
For all positive integers $n$ and $y$, we have
$$
F_{n+5}+F_n\ne y^5.
$$
\end{theorem}

The paper is organized as follows. In Section~2, we establish the
arithmetic facts needed for the two quadratic fields. In Sections~3
and~4, we solve the Diophantine equations in \eqref{eq:pair}.
Section~5 gives the application to the Fibonacci sequence. We then
discuss the scope of the method. The Magma and PARI/GP scripts and
their outputs are given in the appendices.

\section{Preliminaries}
\begin{lemma}\label{lem:parity}
Suppose $Z>0$ and $X^2-5Z^5=44$ or $X^2+44=5Z^5$. Then $X$ and $Z$ are odd, $11\nmid XZ$, and the conjugate ideals generated by $X\pm2\sqrt{11}$ (in the first case) or $X\pm2\sqrt{-11}$ (in the second case) are coprime.
\end{lemma}
\begin{proof}
If $Z$ is even, then the first equation gives $X^2\equiv12\pmod{32}$, and the second gives $X^2\equiv20\pmod{32}$. Neither is possible. Hence $Z$ and $X$ are odd. If $11\mid X$, either equation modulo $11$ gives $11\mid Z$, whereas the left side $X^2\mp44$ has $11$-adic valuation $1$ and $5Z^5$ has valuation at least $5$. Thus $11\nmid XZ$. Suppose a prime ideal $\mathfrak p$ divides both conjugate factors. Then $\mathfrak p$ divides their difference $4\sqrt{\pm11}$, so the rational prime below $\mathfrak p$ is either $2$ or $11$. If $\mathfrak p$ lies above $2$, the congruence $X+2\sqrt{\pm11}\equiv X\pmod{\mathfrak p}$ contradicts the fact that $X$ is odd. If $\mathfrak p$ lies above $11$, then $\sqrt{\pm11}\in\mathfrak p$ and hence $X\in\mathfrak p$, contradicting $11\nmid X$. Thus the two ideals are coprime.
\end{proof}
\section{Classification of $X^2-5Z^5=44$}
Let $K_+=\mathbb Q(r)$, $r^2=11$. Its ring of integers is $\mathbb Z[r]$, with discriminant $44$. Its class number is one: the Minkowski bound is $\sqrt{44}/2<4$, the prime $3$ is inert since $11$ is a nonsquare modulo $3$, and the unique prime ideal above the ramified prime $2$ is principal since $N(3+r)=-2$. Hence every ideal class has a principal representative. The unit group is $\{\pm\varepsilon^m:m\in\mathbb Z\}$, where $\varepsilon=10+3r$. Indeed, $x^2-11y^2=-1$ is impossible modulo $11$, and $10+3r$ is the least unit greater than $1$: for a unit $x+yr>1$ of norm $1$, its conjugate equals its reciprocal, hence $x=(u+u^{-1})/2>0$ and $y=(u-u^{-1})/(2r)>0$. The equation $x^2-11y^2=1$ has no integral solutions for $y=1,2$, and $y\ge3$ implies $x+yr\ge10+3r$. Also $5=(4+r)(4-r)$.
\begin{lemma}\label{lem:plus}
Every integer solution of $X^2-5Z^5=44$ with $Z>0$, up to replacing $X$ by $-X$, gives a solution of one of five quintic Thue equations
$$d_k A(a,b)+c_k B(a,b)=2,\qquad k\in\{-2,-1,0,1,2\},$$
where
\begin{align*}
A(a,b)&=a^5+110a^3b^2+605ab^4,\\
B(a,b)&=5a^4b+110a^2b^3+121b^5,
\end{align*}
and the pairs $(c_k,d_k)$ are
$$\begin{array}{c|rrrrr}k&-2&-1&0&1&2\\\hline c_k&136&7&4&73&1456\\d_k&-41&-2&1&22&439\end{array}.$$
Conversely a solution $(a,b)$ gives $Z=a^2-11b^2$ and $X=c_k A(a,b)+11d_k B(a,b)$; it gives a positive-$Z$ solution whenever $Z>0$.
\end{lemma}
\begin{proof}
The factors of $(X+2r)(X-2r)=5Z^5$ generate coprime ideals by Lemma~\ref{lem:parity}. As $5$ splits into the two distinct principal prime ideals $(4+r)$ and $(4-r)$, unique ideal factorization and class number one yield
$$X+2r=u(4+\delta r)(a+br)^5,\quad \delta\in\{\pm1\},\quad u\in\mathcal O_{K_+}^{\times}.$$
The sign of $u$ can be absorbed into the fifth power, and $\varepsilon^5$ can be absorbed into $(a+br)^5$. Thus $u=\varepsilon^k$, $k\pmod5$. Taking conjugates and negating the equation changes $X$ to $-X$ and sends $(k,\delta)$ to $(-k,-\delta)$, so $\delta=1$ suffices. Use $k=-2,-1,0,1,2$ as representatives. Expand $(a+br)^5=A+rB$ and $\varepsilon^k(4+r)=c_k+d_kr$. Comparing coefficients of $r$ gives the Thue equations, and taking norms gives $Z=a^2-11b^2$.
\end{proof}
\begin{theorem}\label{thm:plus}
The only integer solutions of $X^2-5Z^5=44$ with $Z>0$ are $(X,Z)=(\pm7,1)$.
\end{theorem}
\begin{proof}
The complete Magma Thue solver (\texttt{Solutions(T,2 : Exact := true)}) for the five equations of Lemma~\ref{lem:plus} returns, in the order $k=-2,-1,0,1,2$, the solution sets
$$\varnothing,\quad\{(-1,0)\},\quad\varnothing,\quad\varnothing,\quad\varnothing.$$
The sole representative gives $Z=1$ and $X=-7$. Its negative gives $X=7$. Both satisfy the equation. Appendix~\ref{app:magma} provides the complete computation.
\end{proof}
\section{Classification of $X^2+44=5Z^5$}
Let $K_-=\mathbb Q(s)$, $s^2=-11$, and $\omega=(1+s)/2$. Then $\mathcal O_{K_-}=\mathbb Z[\omega]$, $\omega^2=\omega-3$, and the discriminant is $-11$. The imaginary quadratic Minkowski bound is $(2/\pi)\sqrt{11}<3$; the only possible nontrivial prime ideal norm below this bound is $2$, but $2$ is inert because $T^2-T+3$ is irreducible modulo $2$. Hence the class number is one. The unit group is $\{\pm1\}$, and $5=\pi\bar\pi$ with $\pi=1+\omega$ and $N(\pi)=5$.
\begin{lemma}\label{lem:minus}
Every solution of $X^2+44=5Z^5$, $Z>0$, up to replacing $X$ by $-X$, gives an integer solution of
$$a^5+10a^4b-10a^3b^2-70a^2b^3-20ab^4+17b^5=4.$$
\end{lemma}
\begin{proof}
By Lemma~\ref{lem:parity}, the ideals generated by $X\pm2s$ are coprime. Principal ideal factorization and $5=\pi\bar\pi$ give
$$X+2s=\pm\pi(a+b\omega)^5\quad\text{or}\quad X+2s=\pm\bar\pi(a+b\omega)^5.$$
Absorb the sign into the fifth power. Conjugation followed by negation interchanges the two cases and sends $X$ to $-X$, so it suffices to use $\pi$. Writing
\begin{align*}
(a+b\omega)^5&=A_-(a,b)+B_-(a,b)\omega,\\
A_-(a,b)&=a^5-30a^3b^2-30a^2b^3+30ab^4+15b^5,\\
B_-(a,b)&=5a^4b+10a^3b^2-20a^2b^3-25ab^4+b^5,
\end{align*}
we obtain $\pi(A_-+B_-\omega)=(A_--3B_-)+(A_-+2B_-)\omega$. Since $X+2s=(X-2)+4\omega$, comparison of coefficients yields $A_-+2B_-=4$, which is the stated Thue equation. Taking norms gives $Z=a^2+ab+3b^2$.
\end{proof}
\begin{theorem}\label{thm:minus}
The equation $X^2+44=5Z^5$ has no integer solutions with $Z>0$.
\end{theorem}
\begin{proof}
Magma's complete Thue solver (\texttt{Solutions(T,4 : Exact := true)}) returns the empty set for the equation in Lemma~\ref{lem:minus}. See Appendix~\ref{app:magma}.
\end{proof}
\section{Application to Fibonacci fifth powers}
Let $F_0=0$, $F_1=1$, and $F_{n+2}=F_{n+1}+F_n$. Put $S_n=F_{n+5}+F_n$. Then $S_0=5$, $S_1=9$, and $S_{n+2}=S_{n+1}+S_n$. The standard determinant identity for a Fibonacci-type recurrence gives
$$S_{n+1}^2-S_{n+1}S_n-S_n^2=11(-1)^n,$$
and hence
$$(2S_{n+1}-S_n)^2-5S_n^2=44(-1)^n.$$
\begin{proof}[Proof of Theorem~\ref{thm:intro}]
If $S_n=y^5$, put $X=2S_{n+1}-S_n$ and $Z=y^2>0$. For even $n$, Theorem~\ref{thm:plus} forces $Z=1$, hence $S_n=1$, contrary to $S_n\ge9$. For odd $n$, Theorem~\ref{thm:minus} yields an immediate contradiction.
\end{proof}
\section{Discussion and scope}\label{sec:discussion}
The fixed gap $5$ permits a particularly simple recurrence-theoretic
reduction.  In fact, $S_n=F_{n+5}+F_n$ itself satisfies the Fibonacci
recurrence and has initial values $(S_0,S_1)=(5,9)$.  The associated
quadratic invariant has constant absolute value $11$, which leads to
the constants $\pm44$ in \eqref{eq:pair}.  It is the small discriminant
of the resulting quadratic fields, together with their class number
one, that makes the ensuing ideal-theoretic descent economical.

The result is deliberately restricted to fifth powers.  For another
exponent $p$, the same recurrence identity would produce an equation
of the shape $X^2-5y^{2p}=\pm44$, but the corresponding descent and
Thue computations would need to be carried out anew.  Likewise, a
different gap changes the initial values and the quadratic invariant,
so the present six Thue equations do not immediately apply.

Thus, the two Diophantine classifications give a complete answer for the fixed gap and exponent considered here. The argument does not provide a classification for arbitrary gaps or exponents.

\appendix
\section{Reproducible Magma computation}\label{app:magma}
The Magma handbook states that \texttt{Solutions(T,m)} returns all integral pairs solving the associated homogeneous Thue equation, with \texttt{Exact := true} by default \cite{MagmaHandbook}. The original computations used the Magma Online Calculator. The script below is a reproducibility script: it explicitly specifies \texttt{Exact := true} and adds assertions for independent consistency checks; it should not be interpreted as a verbatim transcript of the original interactive session. The displayed solution sets were reported from the original runs. The computations were performed using the Magma Online Calculator; the version and remote server configuration were not recorded. The checks of norms and reconstructed Diophantine equations are algebraic consistency checks, not independent certificates of completeness.
\begin{lstlisting}[basicstyle=\ttfamily\footnotesize,breaklines=true]
ZZ := Integers();
P<x> := PolynomialRing(ZZ);
A := x^5 + 110*x^3 + 605*x;
B := 5*x^4 + 110*x^2 + 121;
ks := [-2,-1,0,1,2];
cs := [136,7,4,73,1456];
ds := [-41,-2,1,22,439];
for i in [1..5] do
    c := cs[i]; d := ds[i];
    assert c^2 - 11*d^2 eq 5;
    f := d*A + c*B;
    assert IsIrreducible(f);
    print "even k =", ks[i], "f =", f;
    T := Thue(f);
    time sols := Solutions(T,2 : Exact := true);
    print "solutions =", sols;
    if i eq 2 then
        assert sols eq [[-1,0]];
    else
        assert #sols eq 0;
    end if;
    for sol in sols do
        a := ZZ!sol[1]; b := ZZ!sol[2];
        aa := a^5+110*a^3*b^2+605*a*b^4;
        bb := 5*a^4*b+110*a^2*b^3+121*b^5;
        assert d*aa+c*bb eq 2;
        z := a^2-11*b^2;
        X := c*aa+11*d*bb;
        assert X^2-5*z^5 eq 44;
        print <a,b,z,X>;
    end for;
end for;

fodd := x^5+10*x^4-10*x^3-70*x^2-20*x+17;
assert IsIrreducible(fodd);
T := Thue(fodd);
time oddsols := Solutions(T,4 : Exact := true);
print "odd solutions =", oddsols;
assert #oddsols eq 0;
for sol in oddsols do
    a := ZZ!sol[1]; b := ZZ!sol[2];
    aa := a^5-30*a^3*b^2-30*a^2*b^3
          +30*a*b^4+15*b^5;
    bb := 5*a^4*b+10*a^3*b^2-20*a^2*b^3
          -25*a*b^4+b^5;
    assert aa+2*bb eq 4;
    X := aa-3*bb+2;
    z := a^2+a*b+3*b^2;
    assert X^2+44 eq 5*z^5;
    print <a,b,z,X>;
end for;
\end{lstlisting}
\noindent Expected solution sets for the even cases are $\varnothing$, $\{(-1,0)\}$, $\varnothing$, $\varnothing$, $\varnothing$, and for the odd case $\varnothing$.
\section{Independent verification with PARI/GP}\label{app:pari}
The six Thue equations were independently evaluated using 64-bit PARI/GP version 2.19.0 (release dated September 28, 2026), on Windows 11 with an Intel Core i7-12800H processor and 32 GB RAM, using \texttt{thueinit(f,1)} and \texttt{thue(T,m)}. The second argument \texttt{1} requests unconditional certification of the number-field data. The actual outputs supplied by the independent run were, for the five real-quadratic cases, respectively,
$$
[],\quad[[-1,0]],\quad[],\quad[],\quad[],
$$
and for the imaginary-quadratic case, \texttt{[]}. These coincide with the Magma outputs.
\begin{lstlisting}[basicstyle=\ttfamily\footnotesize,breaklines=true]
A = x^5 + 110*x^3 + 605*x;
B = 5*x^4 + 110*x^2 + 121;
print("k=-2: ",thue(thueinit(-41*A+136*B,1),2));
print("k=-1: ",thue(thueinit(-2*A+7*B,1),2));
print("k=0: ",thue(thueinit(A+4*B,1),2));
print("k=1: ",thue(thueinit(22*A+73*B,1),2));
print("k=2: ",thue(thueinit(439*A+1456*B,1),2));
fodd = x^5+10*x^4-10*x^3-70*x^2-20*x+17;
print("Odd case: ",thue(thueinit(fodd,1),4));
\end{lstlisting}
\section*{Computational environment and certification}
All six quintic Thue equations were solved using the \emph{Magma Online Calculator}. As Magma ran on a remote server, no local CPU or memory specification is attributed to that computation; its server version and configuration were not recorded. The calculations were independently repeated with \emph{64-bit PARI/GP 2.19.0} (September 28, 2026) on \emph{Windows 11}, using an \emph{Intel Core i7-12800H} processor and \emph{32 GB RAM}. PARI/GP was invoked with \texttt{thueinit(f,1)} for unconditional certification, followed by \texttt{thue(T,m)}. The two systems produced identical solution sets in every case. The classification depends on exhaustive Thue solvers, not a bounded search; the scripts and recorded outputs appear above. The Magma assertions additionally check the reconstructed Diophantine identities, although such checks alone do not establish completeness. For published uses of complete Thue computations in related Diophantine settings, see \cite{McDaniel2011,BMS2006II}. For wider background on effective determination of integral points, see \cite{BugeaudEtAl2008}.

\end{document}